# Generalized Functions and Infinitesimals

J.F.Colombeau, (jf.colombeau@wanadoo.fr)

Abstract: It has been widely believed for half a century that there will never exist a nonlinear theory of generalized functions, in any mathematical context. The aim of this text is to show the converse is the case and invite the reader to participate in the debate and to examine the consequences at an unexpectedly elementary level. The paradox appears as another instance of the historical controversy on the existence of infinitesimals in mathematics, which provides a connection with nonstandard analysis. This text is the written version of a talk at the meeting of Nonstandard Analysis, Pisa, june 2006.

**1-Prerequisites.** If f, g are functions of class $\mathbf{C}^1$ on $\mathbb{R}$ let us recall the formula for integration by parts

$$\int_a^b f'(x)g(x)dx = -\int_a^b f(x)g'(x)dx + f(b)g(b) - f(a)g(a).$$

If φ is a $\mathbf{C}^\infty$ function on $\mathbb{R}$ which vanishes outside a bounded set ("φ(-∞)=0= φ(+∞)") the formula becomes

$$\int f'(x)\varphi(x)dx = -\int f(x)\varphi'(x)dx.$$

The concept of "distribution" [16,29,27] consists of interpreting f' as the linear map

$$\varphi \to -\int f(x)\varphi'(x)dx$$

which makes sense even if f is not differentiable. If Ω is an open set in $\mathbb{R}^N$ we denote by $\mathbf{D'}(\Omega)$ the vector space of all distributions on Ω; we will not need to enter any more into the distributions but simply note that elements of $\mathbf{D'}(\Omega)$ have many of the properties of $\mathbf{C}^\infty$ functions on Ω concerning differentiation.

In short: the concept of distributions allows one to freely differentiate rather irregular functions (that are not differentiable in the classical sense) at the price that their partial derivatives are objects (distributions) that are not usual functions. The typical example is given by the Heaviside function H. This is defined by: H(x) = 0 if x < 0, H(x) = 1 if x > 0 (and H(0) unspecified). H is not differentiable at x = 0 in the classical sense because of the discontinuity there. However its derivative, in the sense of distributions, is the "Dirac delta function" δ: δ(x)=0 if x≠0, δ(0) "infinite" so that

$$\int \delta(x)dx = [H]_{-\infty}^{+\infty} = 1-0 =1.$$

The above explains calculations in which physicists differentiate functions that, like H, cannot be differentiated in the classical sense. Physicists not only differentiate irregular functions, they also mix differentiation and multiplication by treating formally irregular functions as if they were $\mathbf{C}^\infty$ functions. But L. Schwartz claimed [27 p10 of the 1966 edition]: "A general multiplication of distributions is impossible in any theory [of generalized functions], possibly different from distribution theory, where there exist a differentiation and a Dirac delta function". More precisely his theorem states [28]: there does not exist an algebra **A** such that:

1) the algebra $\mathbf{C}^0(\mathbb{R})$ (of all continuous functions on $\mathbb{R}$) is a subalgebra of $\mathbf{A}$ and the function $x \to 1$ is the unit element in $\mathbf{A}$.
2) there exists a linear map $D: \mathbf{A} \to \mathbf{A}$, ("differentiation") such that
   i) D reduces to the usual differentiation on $\mathbf{C}^1$ functions
   ii) $D(uv) = Du.v + u.Dv \quad \forall u,v \in \mathbf{A}$
   iii) $D.D(x \to |x|) \neq 0$ (in practice it has to be $2\delta$).

But later L.Schwartz supported publication of the conflicting viewpoint below [6].

**2-Nonlinear generalized functions** (in its simplest formulation this theory requires only elementary calculus [8 p163, 9 p261, 13 p10, 19, 30, 32 p99]). 25 years ago [6,7] I found a differential algebra $\mathbf{G}(\Omega)$ (i.e. an algebra with internal partial derivatives: ($\partial/\partial x_i \mathbf{G}(\Omega)) \subset \mathbf{G}(\Omega)$) in the situation $\mathbf{C}^\infty(\Omega) \subset \mathbf{D'}(\Omega) \subset \mathbf{G}(\Omega)$ in which the inclusion $\mathbf{D'}(\Omega) \subset \mathbf{G}(\Omega)$ is canonical (i.e. free from arbitrary choices) and

- the partial derivatives $\partial/\partial x_i$ in $\mathbf{G}(\Omega)$ induce those in $\mathbf{D'}(\Omega)$
- the multiplication in $\mathbf{G}(\Omega)$ induces on $\mathbf{C}^\infty(\Omega)$ the usual multiplication of $\mathbf{C}^\infty$ functions.

Starting from the classical differential algebra $\mathbf{C}^\infty(\Omega)$, slightly different variants of $\mathbf{G}(\Omega)$ are obtained by following a pattern similar to that used in the construction of the real numbers by the method of Cauchy sequences of rational numbers. The pattern is as follows:

A = an algebra of appropriate families $(f_i)_{i \in I}$ of $\mathbf{C}^\infty$ functions on $\Omega$;
I = an ideal of A made of those families $(f_i)_{i \in I}$ "close to zero" as "$i \to \infty$" in the index set I;
$\mathbf{G}(\Omega) = A/I$ (= quotient of the algebra A by the ideal I).

The objects in $\mathbf{G}(\Omega)$ can be treated as $\mathbf{C}^\infty$ functions on $\Omega$ (but not always exactly like $\mathbf{C}^\infty$ functions, which explains various inconsistencies encountered by physicists from "formal" calculations).

The above result seems to be inconsistent with the Schwartz impossibility result, so that at least one of the assumptions in Schwartz's theorem should not hold. The assumption that does not hold is "$\mathbf{C}^0(\Omega)$ is a subalgebra of $\mathbf{G}(\Omega)$" (although $\mathbf{C}^\infty(\Omega)$ is a subalgebra of $\mathbf{G}(\Omega)$). Let f,g be two continuous functions on $\Omega$; one has two products: the classical one f.g and a new one (in $\mathbf{G}(\Omega)$) denoted by f⊙g, which in general are different elements of $\mathbf{G}(\Omega)$: f.g ≠ f⊙g. But they are not so much different since $\forall \varphi \in \mathbf{C}_c^\infty(\Omega)$ (i.e. infinitely differentiable with compact support) the integral $|\int (f.g - f \odot g)(x) \varphi(x)dx|$ (which makes sense naturally in the $\mathbf{G}$-context) is a "generalized real number", nonzero but less than r for any real number r>0: in short it is a <u>nonzero infinitesimal real number: infinitesimal numbers appear here: they were not invited, not welcome but imposed !</u>

The $\mathbf{G}$-theory shows a perfect coherence with classical mathematics thanks to these infinitesimals: if in $\mathbf{G}$ you exclude new objects such as $\delta^2, \delta^3, ...$ and if you identify $G_1, G_2 \in \mathbf{G}(\Omega)$ if $\forall \varphi \in \mathbf{C}_c^\infty(\Omega)$ $\int (G_1-G_2)(x) \varphi(x)dx$ is infinitesimal, then you obtain $\mathbf{D'}(\Omega)$, but you have lost the structure of an algebra ($\mathbf{D'}(\Omega)$ is only a vector space). <u>Therefore the Schwartz non-existence result is based on the refusal of infinitesimals.</u>

Let us show with the example of the Heaviside function why the classical product has to be - infinitesimally-changed ( also for $\mathbf{C}^p$, p finite, but not for $\mathbf{C}^\infty$ functions).

→Compute the integral
$$I=\int (H^2(x)-H(x))\cdot H'(x)dx$$
where H denotes the Heaviside step function and H' its derivative (the Dirac delta distribution). H may be considered as an idealization of a $\mathbf{C}^1$ function with a jump from the value 0 to the value 1 in a tiny interval around x = 0. Thus classical calculations are justified: $I = [H^3/3 - H^2/2]_{-\infty}^{+\infty} = 1/3 - 1/2 = -1/6$. This implies that $H^2 \neq H$ (since $I \neq 0$): $H^2$ and H differ at x = 0, precisely where H' takes an "infinite value", and this undefined form $0 \times \infty$ gives here the value -1/6 after integration. Therefore the classical formula $H^2 = H$ has to be considered as erroneous in a context suitable to compute I. But it holds in the sense that $\forall \varphi \in \mathbf{C}_c^\infty(\Omega)$ $\int (H^2(x)-H(x))\cdot\varphi(x)dx$ is infinitesimal. We denote this as a "<u>weak equality</u>" or "<u>association</u>" $H^2 \approx H$ (notice that $H^2 \neq H$). In $\mathbf{G}(\Omega)$ we say that $G_1 \approx G_2$ ("$G_1$ <u>associated to</u> $G_2$") if $\forall \varphi \in \mathbf{C}_c^\infty(\Omega)$ $\int (G_1-G_2)(x)\,\varphi(x)dx$ is infinitesimal. If $T_1, T_2$ are two distributions on $\Omega$ one proves that $T_1 \approx T_2$ in $\mathbf{G}(\Omega)$ implies $T_1 = T_2$ in $\mathbf{D}'(\Omega)$. For physical applications it will be basic to have in mind that although there is only one Heaviside distribution, there is an infinity of Heaviside like objects in $\mathbf{G}$ (in particular any $H^N$); all of them are called "Heaviside generalized functions"; of course the same holds for their derivatives: "Dirac delta generalized functions" [5 p49, p52, 8 p32, p47].

→Different Heaviside functions are also very concretely imposed by physics: even at an obvious qualitative level the depiction of an elasto-plastic shock wave requires very different Heaviside functions for different physical variables [5 p120, 8 p101, p106].

→Here is a simplified version of Schwartz's proof in which we assume that the algebra of step functions is a subalgebra of $\mathbf{A}$ (instead of $\mathbf{C}^0(\mathbb{R})$, so as to permit a much simpler proof). In the classical algebra of step functions, therefore in $\mathbf{A}$ from our specific assumption: $H^2 = H$ and $H^3 = H$. By differentiation (we note $DH = H'$): $2HH' = H'$ and $3H^2H' = H'$. Since $H^2 = H$ in the algebra $\mathbf{A}$, $3H^2H' = 3HH'$, thus $3HH' = H'$. Thus we have at the same time that $HH' = 1/2 \cdot H'$ and $HH' = 1/3 \cdot H'$, which implies $H' = 0$. ∎
In $\mathbf{G}$, one has $H^2 \approx H$ and $H^3 \approx H$ which gives by differentiation $2HH' \approx H'$ and $3H^2H' \approx H'$ but $H^2H'$ is not weakly equal to $HH'$ since in general weak equality $\approx$ is not preserved by multiplication, which fortunately stops the above calculation before its end.

Therefore if we denote by $\mathbf{D}'(\Omega)^\approx$ the vector space of those elements of $\mathbf{G}(\Omega)$ associated to a distribution one has the set of inclusions:
$$\mathbf{C}^\infty(\Omega) \subset \mathbf{D}'(\Omega) \subset \mathbf{D}'(\Omega)^\approx \subset \mathbf{G}(\Omega).$$
The partial derivatives $\partial/\partial x_i$ are internal in all these four spaces, but only $\mathbf{C}^\infty(\Omega)$ and $\mathbf{G}(\Omega)$ are algebras. If N>1 $H^N$ is in $\mathbf{D}'(\Omega)^\approx$ and not in $\mathbf{D}'(\Omega)$, $\delta^N$ is in $\mathbf{G}(\Omega)$ and not in $\mathbf{D}'(\Omega)^\approx$.

**3-Infinitesimals in continuum mechanics** (main references: [8,21]). About 1980 there appeared projectiles that destroyed in one shot all existing models of battle tanks, thus implying the need to design new armour. Impacts last only a few microseconds: it is impossible to perform detailed experiments and therefore numerical simulations were indispensable in the design of the new armour requested. "Trivial discretizations", such as $f'(x) \# (f(x+h)-f(x))/h$, need to be complemented by "artificial viscosity" which in turn makes the numerical schemes degenerate too quickly. This shows the need to use good quality schemes – here Godunov schemes – for the system of solid mechanics. Godunov schemes are based on explicit solutions of the particular case of the Cauchy problem when the initial data $x \mapsto u(x,0)$ is constant on both sides of a discontinuity (this is called the Riemann problem).

The equations one needs to solve are the equations of continuum mechanics for solids, see [8 p14]. A very simplified model used to demonstrate the method is given by the system of 3 equations that physicists state as

$$\rho_t + (\rho u)_x = 0, \quad (\rho u)_t + (\rho u^2)_x = \tau_x, \quad \tau_t + u\tau_x = u_x,$$

where $\rho = \rho(x,t)$=volumic mass, $u=u(x,t)$=velocity, $\tau=\tau(x,t)$=stress and indices t,x denote partial derivatives with respect to t and x. The first equation is the equation of mass conservation, the second one is the equation of momentum conservation, and the third one is a state law of elastic solids in fast deformation (with coefficient chosen equal to 1). Because of the presence of the term $u\tau_x$ in the third equation there appears a product of the kind H.δ when one seeks a solution of the Riemann problem. The product H.δ does not make sense within distributions, so one therefore works in the **G**-setting. However, as already noted, in the **G**-setting one has to consider an infinity of Heaviside like objects: the formula

$$w(x,t) = w_l + (w_r - w_l)H_w(x - ct),$$

(where $H_w$ is some Heaviside generalized function), expresses the fact that the physical variable w takes the value $w_l$ if x<ct and the value $w_r$ if x>ct (discontinuity travelling at constant speed c). From the above we see that the product $u\tau_x$ involves the product $H_u.(H_\tau)'$ that can take different values depending on $H_u$ and $H_\tau$. The product $u\tau_x$ is ambiguous in the **G**-context, in absence of further information. There are three natural ways of dealing with this ambiguity.

<u>First idea</u>: state all 3 equations of the simplified model with the algebraic equality (=) in **G**. There does not exist a discontinuous solution [8 §332 and p72]; this is not acceptable!

<u>Second idea</u>: state all 3 equations with the weak equality ≈ in **G**, i.e.:

$$\rho_t + (\rho u)_x \approx 0, \quad (\rho u)_t + (\rho u^2)_x \approx \tau_x, \quad \tau_t + u\tau_x \approx u_x.$$

There are infinitely many different solutions for a given initial condition [8 p69, 9 p265]; this is not acceptable!

<u>Third idea</u>: choose an intermediate statement on physical ground. Physicists have observed that shock waves have an "infinitesimal width" of the order of magnitude of a few hundred crystalline meshes. One can therefore isolate (by thought) small volumes inside this width where conservation laws apply. This suggests to state them with (algebraic) = in **G**. On the other hand the state law (3rd equation) has been checked only on a material at rest, i.e. on both sides of the shock wave: by analogy with H² ≈ H this suggests that one should state it with the weak equality ≈ in **G**. Thus one states in **G** the system of 3 equations as:

$$\rho_t + (\rho u)_x = 0, \quad (\rho u)_t + (\rho u^2)_x = \tau_x, \quad \tau_t + u\tau_x \approx u_x.$$

With this formulation of the problem one obtains the desired existence-uniqueness result and explicit formulas [8 p72]. This method has been used in [21] for numerical simulations of collisions. For contact discontinuities this physical idea suggests to state the state laws with = in **G** since there is no fast deformation: one obtains the desired existence-uniqueness result.

**4-Infinitesimals in general relativity** (main references: [13,32]). For charged ultrarelativistic (moving at lightspeed) black holes physicists have noticed the rather unexpected result that the electromagnetic field vanishes but its energy-momentum tensor does not. This "absurdity" should be clarified: has it its origin in a mathematical mistake or does it point out a breakdown of physics? Rigorous calculations in the **G**-setting ([30]) show that the field looks like $\sqrt{\delta}$, the square root of a Dirac delta generalized function, while its energy momentum tensor (involving the square of the field) is δ-like. Since $\sqrt{\delta} \approx 0$ the field is infinitesimal but nonzero, thus permitting its energy momentum tensor to be nonvanishing. Then the physically unsatisfactory situation of a vanishing field with nonzero energy momentum tensor is mathematically perfectly clarified by the **G**-setting: the paradox was due to a lack of rigor in the mathematical calculations.

In suitable coordinates impulsive gravitational waves can be represented as follows [13 p435, 32 p104]: the space time is flat except for a hypersurface u=0 where a $\delta$-like impulse modelling a gravitational shock wave is located. They are treated by the "scissors and paste" method of Penrose: space time is divided into two halves by removal of the hypersurface u=0, then the two halves are joined together by the "Penrose junction conditions". The corresponding mathematical calculations involve nonlinear generalized functions (ill-defined products of distributions within distribution theory) due to the nonlinearity of Einstein equations and the presence of a Dirac delta function in the space-time metric. It is shown in [20] that these calculations represent a well defined diffeomorphism in the **G**-sense.

Einstein's equations are nonlinear and classically require the $C^2$-setting. However, one frequently encounters less regular "solutions" (discontinuous, $\delta$-like,…). Geroch and Traschen [12] have shown that a physically sensible and mathematically sound setting strictly within distribution theory excludes the description of gravitational sources concentrated on submanifolds of codimension greater than 1 in space-time: physically interesting sources like cosmic strings or point particles are strictly excluded. This fact and the examples above show the need for a semi-Riemannian geometry based on nonlinear generalized functions [32 §8]. It is currently applied in the context of the singularity theorems of general relativity with the hope of distinguishing between "weak singularities", which correspond to physically realistic models which can be mathematically tackled within the **G**-setting and "real singularities" where the space time geometry breaks down and new physics is needed [32 p110].

**5-Questions.** The nonlinear theory of generalized functions is inevitably based on infinitesimals, it applies to physics and of course also to mathematics motivated by physics (may be more than 500 articles on applications). Therefore this theory contributes to the debate on existence of infinitesimals in mathematics [1,18,23,24,26,33].

1. Should it be considered as some kind of nonstandard analysis? See [16 chap6 and references there, 25 p258, 34]: infinitesimals, compatibility of both theories.
2. What improvements could nonlinear generalized functions bring to nonstandard analysis? Maybe an enlargement of its field of applications [2,11].
3. **What improvements could nonstandard analysis bring to nonlinear generalized functions? Maybe tools from logic such as ultrapowers, transfer and saturation principles [16,18,23,33].**

**6-A definition.** We give a definition of a simplified algebra $\mathbf{G}_S(\Omega)$ (the inclusion of $\mathbf{D}'(\Omega)$ into $\mathbf{G}_S(\Omega)$ is not canonical). I=]0,1]; $\mathbf{A}=\{(u_\varepsilon)_{\varepsilon \in I}, u_\varepsilon \in \mathbf{C}^\infty(\Omega)$, such that $\forall K \subset\subset \Omega$ $\forall \alpha \in IN^n, \exists N \in IN$ with $\sup_{x \in K} |\partial^\alpha u_\varepsilon(x)| = O(\varepsilon^{-N})$ as $\varepsilon \to 0$ }; $\mathbf{I}=\{(u_\varepsilon)_{\varepsilon \in I} \in \mathbf{A}$, such that $\forall K \subset\subset \Omega \ \forall q \in IN \ \sup_{x \in K} |u_\varepsilon(x)| = O(\varepsilon^q)$ as $\varepsilon \to 0$ }; $\mathbf{G}_S(\Omega) = \mathbf{A}/\mathbf{I}$; $\subset\subset$ means compact.

From a nice remark due to Grosser [13 p11] we do not need to introduce $\partial^\alpha$ in the definition of **I** as done in most texts. The inclusion $\mathbf{C}^\infty(\Omega) \subset \mathbf{G}_S(\Omega)$ is obtained by choosing $u_\varepsilon = f$ $\forall \varepsilon$ if $f \in \mathbf{C}^\infty(\Omega)$. No need to know anything on distributions. Is this definition really more complicated than the one of the real numbers? It reflects the opinion that "generalized functions" are "idealizations" defined from classical functions (the $u_\varepsilon$'s) by some kind of "abstract approximation process" allowing good mathematical properties: a product of conceptual thought done on purpose to be very applicable to physical reality (like the real numbers). See [32 p111] for the respective roles of $\mathbf{G}_S(\Omega)$ and $\mathbf{G}(\Omega)$ in physics ("special" and "full" algebra respectively in the terminology of [13,32]).

**7-A "game"**. In a classical physical situation (no need to enter into general relativity or quantum mechanics) introduce an irregularity: consider a medium made of different layers (=use of a step function to model it), or consider a population concentrated at a point (= use of a Dirac delta function to model it). The models often show products of the kind H.δ [4,17,25 p161-164] or $\delta^2$ [10,25 p170-180]. Anybody can find examples and solve them: minimization of the integrals (brachistochron and catenoid)

$$B = \int_0^1 \chi(x)((1+u'(x)^2)/u(x))^{1/2} dx \quad , \quad C = \int_0^1 \chi(x) \cdot u(x) \cdot (1+u'(x)^2)^{1/2} dx$$

leads to ODEs (Euler-Lagrange) involving products of distributions if χ is a discontinuous weight function: χ might model a danger or a price; consider also χ(x)=1+ δ(x-1/2). One can do immediately new mathematics (and deeper formulations of physics as in §3) that work amazingly well. Whether models in the **G**-context deserve a detailed study depends on their use in <u>engineering</u>: §3, hydrodynamics [3,4,8,17,21], hurricanes [22], earthquakes [14,15], or in <u>physics</u>: §4, Schwarzschild and Kerr spacetimes, ultrarelastivistic black holes, geodesics for irregular spacetimes, cosmic strings, "generalized hyperbolicity", "physically reasonable" and "true gravitational" singularities [13,19,20,30,31,32 and references there].

The A. is indebted to G. Berger, R. Steinbauer and J.Vickers for help in the preparation of this paper.
jf.colombeau@wanadoo.fr ; 33 rue de la Noyera, pavillon 17, 38090, Villefontaine, France.
**References.**
[1]L.Arkeryd. Nonstandard Analysis. Amer.Math.Monthly 112,2005,p926-928.
[2]L.Arkeryd,N.Cutland,CW.Henson(eds).Nonstandard Analysis,Theory and Applications. Springer, 1997.
[3]R.Baraille,G.Bourdin,F.Dubois,A.Y.LeRoux. Une version à pas fractionnaires du schéma de Godunov pour l'hydrodynamique.Comptes Rendus Acad.Sci.Paris 314,1,1992,p147-152.
[4]F.Berger,JF.Colombeau. Numerical solutions of one –pressure models in multifluid flows. SIAM J. Numer. Anal.32,4,p1139-1154,1995.
[5]H.A.Biagioni. A Nonlinear Theory of Generalized Functions. Lecture Notes in Math 1421, Springer Verlag, 1990.
[6]J.F.Colombeau. A general multiplication of distributions. Comptes Rendus Acad.Sci.Paris 296,1983,p357-360, and subsequent notes presented by L.Schwartz.
[7]J.F.Colombeau. A multiplication of distributions, J.Math.Ana.Appl. 94,1,1983, p96-115. New Generalized Functions and Multiplication of Distributions. North Holland, Amsterdam, 1984. Elementary Introduction to New Generalized Functions. North Holland, Amsterdam, 1985.
[8]J.F.Colombeau. Multiplication of Distributions. Lecture Notes in Maths 1532. Springer Verlag 1992. See also:The elastoplastic shock problem as an example of the resolution of ambiguities in the multiplication of distributions J. Math Phys 30, 90, 1989, p2273-2279.
[9]J.F.Colombeau. Multiplication of distributions. Bull. AMS. 23,2, 1990, p251-268.
[10]J.F.Colombeau,M.Oberguggenberger. Hyperbolic system with a compatible quadratic term, delta waves and multiplication of distributions. Comm.in PDEs 15,7,1990, p905-938.
[11]N.Cutland (ed.). Nonstandard Analysis and its Applications.Cambridge Univ.Press, 1988.
[12]R.Geroch, J.Traschen. Strings and other distributional sources in general relativity. Phys. Rev D, 36 (4), p1017-1031, 1987.

[13] M.Grosser, M.Kunzinger, M.Oberguggenberger, R.Steinbauer. Geometric Theory of Generalized Functions with Applications to General Relativity. Kluwer, Dordrecht-Boston-New York, 2001.

[14] G.Hoermann, M.de Hoop. Geophysical modelling with Colombeau functions. arXiv AP/0104007.

[15] G.Hoermann, M.de Hoop. Detection of wave front set perturbations via correlations: foundations for wave equation tomography. Appl.Anal.81,2002,p1443-1465.arXiv math-ph/0104003.

[16] R.F.Hoskins, J.Sousa-Pinto. Distributions, Ultradistributions and Other Generalised Functions. Ellis-Horwood, New York-London, 1994, and references there (Todorov,..)

[17] J.Hu. The Riemann problem for pressureless fluid dynamics in the Colombeau algebra. Comm. Math. Phys. 194,1, p191-205,1998.

[18] AE.Hurd, PA.Loeb. An Introduction to Nonstandard Real Analysis. Academic Press,1985.

[19] M.Kunzinger, R.Steinbauer. A rigourous solution concept for geodesics and geodesic deviation equations in impulsive gravitational waves. J. Math. Phys. 40, p1479-1489,1999. arXiv gr-qc/9806009.

[20] M.Kunzinger, R Steinbauer. A note on the Penrose junction conditions. Class. Quant. Grav. 16, p1255-1264, 1999. arXiv gr-qc/9811007.

[21] A.Y.LeRoux. Phd thesises held at the University of Bordeaux 1, to be found from the university library: Cauret 1986, Adamczewski 1986, De Luca 1989, Arnaud 1990, Baraille 1992.

[22] A.Y.LeRoux, M.N.LeRoux, J.A.Marti. A mathematical model for hurricanes. Comptes Rendus Acad.Sci.Paris 339,1, p313-316,2004.

[23] W.A.J.Luxemburg. What is Nonstandard Analysis? Amer. Math. Monthly, supplement, p38-67.

[24] F.A.Medvedev. Nonstandard Analysis and the history of classical analysis. Amer .Math. Monthly. 1998, p659-664.

[25] M.Oberguggenberger. Multiplication of Distributions and Applications to PDEs. Pitman Research Notes in Math 259, 1992, Longman.

[26] A.Robinson. Nonstandard Analysis. North-Holland 1966. Princeton Univ.Press. 1996.

[27] L.Schwartz .Théorie des Distributions, Hermann, Paris, numerous editions.

[28] L.Schwartz. Sur l'impossibilité de la multiplication des distributions. Comptes Rendus Acad. Sci. Paris 239, 1954, p847-848.

[29] S.L.Sobolev. Méthode nouvelle à résoudre le problème de Cauchy. Math. Sbornik 1,1936, p39-71.

[30] R.Steinbauer. The ultrarelastivistic Reissner-Nordström field in the Colombeau algebra. J. Math. Phys. 38, p1614-1622, 1997. arXiv gr-qc/9606059.

[31] R.Steinbauer. Geodesic and geodesic deviation for impulsive gravitational waves. J. Math. Phys. 39, p2201-2212, 1998. arXiv gr-qc/9710119

[32] R.Steinbauer, JA.Vickers. The use of generalized functions and distributions in general relativity. Class.Quant.Grav.23, pR91-114, 2006, numerous references there, many in arXiv. .

[33] KD.Stroyan, WAJ.Luxemburg. Introduction to the Theory of Infinitesimals. Academic Press,1976.

[34] T.D. Todorov. Colombeau's generalized functions … in Stankovic, Pap, Pilipovic, Vladimirov editors. Generalized Functions, Convergence Structures, and their Applications, Plenum Press, New York 1988, p327-339.